\documentclass[12pt]{article}
\usepackage{graphicx}
\usepackage{amssymb}
\newtheorem{theorem}{Theorem}[section]

\newtheorem{proposition}[theorem]{Proposition}

\begin{document}
 \def\g{\mbox{${\bf g}$}}
 \def\F{\mbox{${\cal F}$}}

\newcounter{rown}
\def\bl{\setcounter{rown}{\value{equation}}
        \stepcounter{rown}\setcounter{equation}0
        \def\theequation{\thesection.\arabic{rown}\alph{equation}}
        }
\def\el{\setcounter{equation}{\value{rown}}
        \def\theequation{\thesection.\arabic{equation}}
        }
\def\sec{\setcounter{equation}0}
\renewcommand{\theequation}{\thesection.\arabic{equation}}

\title{
Chains of extended Jordanian twists for \\
Lie superalgebras 
\thanks{This work was supported by Russian
Foundation for Fundamental Research, grant No. RFBR-02-01-00668,
and CRDF RMI-2334-MD-02.}
\thanks{The talk given on The International
Workshop "Supersymmetries and Quantum Symmetries -
SQS'03",(Russia, Dubna, July 24-31, 2003).}}

\author{V.N. Tolstoy\\
Institute of Nuclear Physics, Moscow State University \\
119992 Moscow \& Russia (e-mail: tolstoy@nucl-th.sinp.msu.ru)}
\date{}
\maketitle

\begin{abstract}
Two type of superization of the Jordanian $r$-matrix for the Lie
algebra $\mathfrak{sl}(2)$ are considered. One type is associated
with the Lie superalgebra $\mathfrak{sl}(1|1)$ and another type is
associated with the orthosymplectic Lie superalgebra
$\mathfrak{osp}(1|2)$. Extended Jordanian $r$-matrices of maximal
order are obtained for the basic complex Lie superalgebras
$\mathfrak{sl}(m|n)$ and $\mathfrak{osp}(M|2n)$, and a general
procedure for construction of corresponding chains of extended
Jordanian twists is given. We also find a relation between the
extended Jordanian twist and automorphism which gives trivial
coproduct for a subalgebra provided the subalgebra is a kernel of
the cobracket for the corresponding $r$-matrix.
\end{abstract}

\section{Introduction}
The Drinfeld's quantum group theory roughly includes two classes
of Hopf algebras: quasitriagular and triangular. The (standard)
$q$-deformation of simple Lie algebras belongs to the first class.
The simplest example of the triangular (non-standard) deformation
is the Jordanian deformation of  $\mathfrak{sl}(2)$ (e.g., see
\cite{O}). In the case of simple Lie algebras of rank $\geq2$ some
non-standard deformations were constructed by  Kulish, Lyakhovsky
et al. \cite{KLM}--\cite{AKL}. These deformations are described by
chains of twists which are an extension of the Jordanian twist.
Full chains of extended Jordanian twists were constructed for all
complex Lie algebras of the classical series $A_n$, $B_n$, $C_n$
and $D_n$.

In this paper a generalization of these results on the supercase
is given. Namely, we consider two type of superization of the
Jordanian $r$-matrix for the Lie algebra $\mathfrak{sl}(2)$. One
type is associated with the Lie superalgebra $\mathfrak{sl}(1|1)$
and another type is associated with the orthosymplectic Lie
superalgebra $\mathfrak{osp}(1|2)$. Extended Jordanian
$r$-matrices of maximal order are obtained for the basic complex
Lie superalgebras $\mathfrak{sl}(m|n)$ and $\mathfrak{osp}(M|2n)$,
and a general procedure for construction of corresponding chains
of the extended Jordanian twists is given. The super-Jordanian
deformation of the Lie superalgebra $\mathfrak{osp}(1|2)$ was
found in \cite{BLT}.  We also find a relation between the extended
Jordanian twist and automorphism which gives trivial coproduct for
a subalgebra provided the subalgebra is a kernel of the cobracket
for the corresponding $r$-matrix.

\setcounter{equation}{0}
\section{Classical $r$-matrices of Jordanian type}

First we consider some notations and definitions concerning classical
$r$-matrices.

Let $\mathfrak{g}$ be any finite-dimensional simple Lie
superalgebra  then
$\mathfrak{g}=\mathfrak{n_-^{}}\oplus\mathfrak{h}\oplus\mathfrak{n_+^{}}$,
where $\mathfrak{n_{\pm}^{}}$ are maximal nilpotent subalgebras
and $\mathfrak{h}$ is a Cartan subalgebra. The subalgebra
$\mathfrak{n_+^{}}$ ($\mathfrak{n_-^{}}$) is generated by the
positive (negative) root vectors $e_{\beta}^{}$ ($e_{-\beta}^{}$ )
for all $\beta\in\triangle_+(\mathfrak{g})$. The symbol
$\mathfrak{b}_{+}^{}$ will denote the Borel subalgebra of
$\mathfrak{g}$,
$\mathfrak{b}_{+}^{}:=\mathfrak{h}\oplus\mathfrak{n}_+^{}$. Let a
Cartan element $h_{\theta}^{}\in \mathfrak{h}$ and a homogeneous
element $e_{\theta}^{}\in \mathfrak{n}_{+}$ satisfies the relation
\begin{equation}
[h_{\theta}^{},\,e_{\theta}^{}]=e_{\theta}^{}\qquad
\bigl(\deg(h_{\theta}^{})=0,\;\; \deg(e_{\theta}^{})=0,\,{\rm
or}\;\,1\bigr)~. \label{jr1}
\end{equation}
Consider the skew-symmetric two-tensor
\begin{equation}
r_{\theta}^{}(\xi)=\xi\,h_{\theta}^{}\wedge e_{\theta}^{}:=
\xi\bigl(h_{\theta}^{}\otimes e_{\theta}^{}-
e_{\theta}^{}\otimes h_{\theta}^{}\bigr)~,\qquad
r_{\theta}^{21}(\xi)=-r_{\theta}^{}(\xi),
\label{jr2}
\end{equation}
where symbol $\xi$ is called a deformation parameter.
We demand that the two-tensor is even, $\deg(r_{\theta}^{}(\xi))=0$.
It means that the deformation parameter $\xi$ also should be homogeneous
and we have
\begin{equation}
\deg(\xi)=\deg(e_{\theta}^{})~,\qquad \xi e_{\theta}^{}=
(-1)^{\deg(\xi)\deg(e_{\theta}^{})}e_{\theta}^{}\xi~.
\label{jr3}
\end{equation}
It easy to check that the two-tensor (\ref{jr2}) with the conditions
(\ref{jr3}) is a classical $r$-matrix, i.e. it satisfies the classical
Yang-Baxter equation (CYBE)
\begin{equation}
[r_{\theta}^{12}(\xi),r_{\theta}^{13}(\xi)+r_{\theta}^{23}(\xi)] +
[r_{\theta}^{13}(\xi),r_{\theta}^{23}(\xi)]=0~.
\label{jr4}
\end{equation}
In the case, when the element $e_{\theta}^{}$ is even,
$\deg(e_{\theta}^{})=0$, and $\xi$ is a complex number, $\xi\in
{\mathbb C}$, the $r$-matrix (\ref{jr2}) is called a Jordanian
classical $r$-matrix and it was intensively discussed in the
literature. In the supercase, when
$\deg(e_{\theta}^{})=\deg(\xi)=1$, the $r$-matrix (\ref{jr2}) will
be also  called a Jordanian classical $r$-matrix.

{\it Remarks.} (i) The condition $\deg(r_{\theta}^{}(\xi))=0$ is
natural since $R=1+r_{\theta}^{}(\xi)+\ldots$, where $R$ is the
universal $R$-matrix. (ii) The second condition (\ref{jr3}) is
consequence of the fact that the co-bracket $\delta(x)$ belongs to
$g\otimes g$, where $\delta(x)$ is defined by
$\xi\delta(x):=[x\otimes1+ 1\otimes x,r_{\theta}^{}(\xi)]$.
(iii)~In the case, when the deformation parameter $\xi$ is odd
("fermion"), there are  two possibilities: $\xi^2=0$ (Grassmann
variable), and $\xi^2\neq0$ (Clifford variable). In this paper we
shall not make this specialization.

{\bf Example 1.} Consider the simplest Lie superalgebra
$\mathfrak{gl}(1|1)$. It is generated by the Cartan-Weyl elements
$h_1:=e_{1-1}$, $h_2:=e_{2-2}$, $e_{1-2}^{}$ and $e_{2-1}^{}$ with
the relations
\begin{equation}
\begin{array}{rclcl}
[h_{i}^{},\,e_{k-l}^{}]\!\!&=\!\!&(\delta_{ik}^{}-\delta_{il}^{})e_{k-l}^{}~,
\quad\{e_{1-2}^{},\,e_{2-1}^{}\}\!\!&=\!\!&h_{1}+h_{2}~, \quad
e_{1-2}^{2}=e_{2-1}^{2}=0~,
\\[7pt]
\deg(h_{1}^{})\!\!&=\!\!&\deg(h_{2}^{})=0~,\qquad\quad
\deg(e_{2-1}^{})\!\!&=\!\!&\deg(e_{1-2}^{})=1~. \label{jr5}
\end{array}
\end{equation}
It is well-known the following classical Drinfeld-Jimbo $r$-matrix for
$\mathfrak{gl}(1|1)$:
\begin{eqnarray}
r_{_{DJ}}^{}(\hbar)&=&\hbar\;(e_{1-2}^{}\otimes
e_{2-1}^{}+e_{2-1}^{}\otimes e_{1-2}^{})\,,
\label{jr6}
\end{eqnarray}
where the parameter $\hbar$ is even, $\deg(\hbar)=0$.

Let $\eta$ be an odd parameter then we can write down one more
additional Jordanian solution:
\begin{eqnarray}
r_{1}^{}(\eta)&=&\eta\;(h_{1}^{}-h_{2}^{})\wedge e_{1-2}^{}~,
\label{jr7}
\end{eqnarray}
where $\deg(\eta)+\deg(e_{1-2})=0$ and $\eta\,e_{1-2}=-e_{1-2}\,\eta$.

{\bf Example 2.} Consider the Lie superalgebra $\mathfrak{osp}(1|2)$.
It is generated by the Cartan-Weyl elements $h$, $v_{\pm}^{}$ and $e_{\pm}^{}$
with the relations
\begin{equation}
\begin{array}{rclcl}
[h,\,v_{\pm}^{}]\!\!&=\!\!&\pm\,\frac{1}{2}\,v_{\pm}\,,\qquad\qquad
\{v_{+}^{},\,v_{-}^{}\}\!\!&=\!\!&-\frac{1}{2}\,h\,,\qquad
e_{\pm}^{}=\pm4\,(v_{\pm}^{})^2\,,\qquad
\\[7pt]
\deg(h)\!\!&=\!\!&\deg(e_{\pm}^{})=0~,\qquad
\deg(v_{\pm}^{})\!\!&=\!\!&1~.
\label{jr8}
\end{array}
\end{equation}
It is well-known the following classical $r$-matrices for
$\mathfrak{osp}(1|2)$ \cite{KR, JS}:
\begin{eqnarray}
r_{_{DJ}}^{}(\hbar)&=&\hbar\;(e_{+}^{}\wedge
e_{-}^{}+2v_{+}^{}\otimes v_{-}^{} +2v_{-}^{}\otimes v_{+}^{})\,,
\label{jr9}
\\[5pt]
r_{1}^{}(\xi)&=&\xi\; h\wedge e_{+}^{}~,\quad
r_{2}^{}(\xi)\;=\;\xi\;(h\wedge e_{+}^{}-2v_{+}^{}\otimes v_{+}^{})\,,
\label{jr10}
\end{eqnarray}
where the parameters $\hbar$ and $\xi$ are even, $\deg(\hbar)=\deg(\xi)=0$.

Let $\eta$ be an odd parameter then we can write down two more additional
solutions:
\begin{eqnarray}
r_{3}^{}(\eta)\!\!&=\!\!&\eta\;h\wedge v_{+}^{}~,\qquad
r_{4}^{}(\eta)\;=\;\eta\;v_{+}^{}\wedge e_{+}^{}~, \label{jr11}
\end{eqnarray}
where $\deg(\eta)+\deg(v_+)=0$ and $\eta\,v_+=-v_+\,\eta$.

Again, let a Cartan element $h_{\theta}^{}\in \mathfrak{h}$ and a
homogeneous root vector $e_{\theta}^{}\in \mathfrak{n}_{+}$ satisfy the
relation (\ref{jr1}). Moreover, let homogeneous elements
$e_{\gamma_{\pm i}}^{}$ indexed by the symbols $i$ and $-i$,
$i\in I=\{1,2,\ldots,N\}$ satisfy the relations
\begin{equation}
\begin{array}{rcl}
[h_{\theta}^{},\,e_{\gamma_i}]\!\!&=&\!\!(1-t_{\gamma_i}^{})\,e_{\gamma_i}~,
\quad[h_{\theta}^{},\,e_{\gamma_{-i}}]\;=\;t_{\gamma_i}^{}\,e_{\gamma_{-i}}~,
\quad (t_{\gamma_i}^{}\in{\mathbb C})~,
\\[7pt]
[e_{\gamma_{\pm
i}},\,e_{\theta}^{}]\!\!&=\!\!&0~,\qquad\qquad\quad\;
[e_{\gamma_k},\,e_{\gamma_{l}}]\;=\;\delta_{k,-l}^{}\,e_{\theta}^{}\,
\quad\bigr(k>l \in I\bigcup\,(-I)\bigl)~,
\\[7pt]
\deg(e_{\theta}^{})\!\!&=\!\!&\deg(e_{\gamma_i})+\deg(e_{\gamma_{-i}})
\quad(\mathop{\rm mod}2)~.
\label{jr12}
\end{array}
\end{equation}
For the Lie superalgebra $\mathfrak{g}$ the brackets
$[\cdot,\,\cdot]$ always denote the super-commutator:
\begin{equation}
[x,y]:=xy - (-1)^{\deg(x)\deg(y)}yx
\label{jr13}
\end{equation}
for any homogeneous elements $x$ and $y$.
Consider the even skew-symmetric two-tensor
\begin{equation}
r_{\theta,N}^{}(\xi)=\xi\,\biggl(h_{\theta}^{}\wedge e_{\theta}^{}+
\sum_{i=1}^{N}(-1)^{\deg(e_{\gamma_i})\,\deg(e_{\gamma_{-i}})}
e_{\gamma_i}^{}\wedge e_{\gamma_{-i}}^{}\biggr)
\label{jr14}
\end{equation}
where
\begin{equation}
\deg(\xi)=\deg(e_{\theta}^{})=\deg(e_{\gamma_i}^{})+\deg(e_{\gamma_{-i}}^{})
\quad(\mathop{\rm mod}2)~.
\label{jr15}
\end{equation}
Moreover we assume that the operation $"\wedge"$ in (\ref{jr14}) is graded:
\begin{equation}
e_{\gamma_i}^{}\wedge e_{\gamma_{-i}^{}}:=e_{\gamma_i}^{}\otimes
e_{\gamma_{-i}}^{}-(-1)^{\deg(e_{\gamma_i})\,\deg(e_{\gamma_{-i}})}
e_{\gamma_{-i}}^{}\otimes e_{\gamma_i}^{}~.
\label{jr16}
\end{equation}
It is not hard to check that the element (\ref{jr14}) satisfies
CYBE and it will be called the extended Jordanian $r$-matrix of
$N$-order. Let $N$ be maximal order, i.e. we assume that another
elements $e_{\gamma_{\pm j}}^{}\in \mathfrak{n}_+^{}$, $j> N$,
which satisfy the relations (\ref{jr12}), do not exist. Such
element (\ref{jr14}) will be called the extended Jordanian
$r$-matrix of maximal order. It is evident that the extended
Jordanian $r$-matrix of maximal order is defined by the elements
$h_{\theta}^{}\in\mathfrak{h}$, $e_{\theta}^{}\in\mathfrak{n}_{+}$
and the Borel subalgebra $\mathfrak{b}_+$. We shall here consider
a special ("canonical") case when $e_{\theta}^{}$ and
$e_{\gamma_{\pm i}}$ ($i=1,2,\cdots,N$) are weight elements with
respect to the Cartan subalgebra $\mathfrak{h}$:
\begin{equation}
[h,\,e_{\theta}^{}]=(h,\theta)\,e_{\theta}^{}\,,\qquad
[h,\,e_{\gamma_{\pm i}}^{}]=(h,\gamma_{\pm i}^{})\,e_{\gamma_{\pm i}}^{}\,
\label{jr17}
\end{equation}
for any $h\in\mathfrak{h}$ and for all $i=1,2,\cdots,N$. Analyzing
the structure of the positive root systems of the complex simple Lie
superalgebras we can see that the very maximal order $N$ of the extended
Jordanian $r$-matrix is associated with the maximal root,
i.e. the root $\theta$ is maximal.

Consider a maximal subalgebra $\mathfrak{b}_+'\in \mathfrak{b}_+$ which
co-commutes with the maximal extended Jordanian $r$-matrix (\ref{jr14}),
$\mathfrak{b}_+':=\mathop{\rm Ker}\delta\in\mathfrak{b_+}$:
\begin{equation}
\xi\delta(x):=[\Delta(x),\,r_{\theta,N}^{}(\xi)]=[x\otimes1+1\otimes x,\,
r_{\theta,N}^{}(\xi)]=0\,\qquad(\forall x\in\mathfrak{b_+'})\,.
\label{jr18}
\end{equation}
Let $r_{\theta_1,N_1}^{}(\xi_1)\in\mathfrak{b}_+'\otimes
\mathfrak{b}_+'$ is also a extended Jordanian $r$-matrix of the form
(\ref{jr14}) with the maximal root $\theta_1\in\mathfrak{h}'$ and
maximal order $N_1$. Then the sum
\begin{equation}
r_{\theta,N;\,\theta_1,N_1}(\xi,\xi_1^{}):=
r_{\theta,N}^{}(\xi)+r_{\theta_1,N_1}^{}(\xi_{1}^{})
\label{jr19}
\end{equation}
is also a classical $r$-matrix, i.e. it satisfies CYBE. This
proposition for superalgebras can be check by direct calculations.
For the case of Lie algebras the proposition was first formulated
in \cite{BD}.

\noindent
Again, we consider a maximal subalgebra $\mathfrak{b}_+^{\prime\prime}\in
\mathfrak{b}_+'$ which co-commutes with the maximal extended Jordanian
$r$-matrix $r_{\theta_1,N_1}^{}(\xi_{1}^{})$ and we construct a extended
Jordanian $r$-matrix of maximal order, $r_{\theta_{2},N_2}^{}(\xi_{2}^{})$.
Continuing this process as result we obtain a canonical chain of subalgebras
\begin{equation}
\mathfrak{b}_+\supset\mathfrak{b}_+'\supset\mathfrak{b}_+^{\prime\prime}\cdots
\supset\mathfrak{b}_+^{(k)} \label{jr20}
\end{equation}
and the resulting $r$-matrix
\begin{equation}
r_{\theta,N;\,\theta_1,N_1;\,\ldots;\theta_k,N_k}^{}
\left(\xi,\xi_{1}^{},\cdots,\xi_{k}^{}\right)=r_{\theta,N}^{}(\xi)+
r_{\theta_1,N_1}^{}(\xi_{1}^{})+\cdots+r_{\theta_k,N_k}^{}(\xi_{k}^{})
\label{jr21}
\end{equation}
is a solution of CYBE. If the chain (\ref{jr20}) is maximal, i.e. it is
constructed in corresponding with the maximal orders $N, N_1,\dots N_k$,
then the $r$-matrix (\ref{jr21}) is called the maximal classical $r$-matrix
of Jordanian type for the Lie algebra $\mathfrak{g}$.

Now we consider examples of maximal classical $r$-matrices of Jordanian
type for the classical Lie superalgebras $\mathfrak{sl}(m|n)$ and
$\mathfrak{osp}(M|2n)$.

{\bf Example 3. Maximal classical r-matrix of Jordanian type
for the Lie superalgebra $A(m|n-1)\simeq\mathfrak{sl}(m|n)$}.
Let $\epsilon_i^{}$ ($i=1,2,\ldots, N:=m+n$) be an orthonormalized
basis of a $N$--dimensional super-Euclidian space $\mathbb R^{(m|n)}$:
$(\epsilon_i^{},\epsilon_j^{})=\pm\delta_{ij}^{}$. In the terms of
$\epsilon_i^{}$ the systems of positive roots, $\Delta_{+}^{}$,
for $\mathfrak{sl}(m|n)$ are presented as follows:
\begin{eqnarray}
\Delta_{+}^{}(\mathfrak{sl}(m|n))\!\!\!\!&=\!\!&\{\epsilon_i^{}-
\epsilon_j^{}\,|\,1\le i<j\le N\}~.
\label{jr22}
\end{eqnarray}
Here some roots are even and another roots are odd. The root
$\epsilon_{1}^{}-\epsilon_{N}^{}$ is maximal. Let us write down
the positive root system $\Delta_{+}^{}(\mathfrak{sl}(m|n))$ in
the following normal ("convex") ordering
\begin{equation}
\begin{array}{rcl}
&(\epsilon_1^{}\!-\epsilon_2^{},\epsilon_1^{}\!-\epsilon_3^{},\ldots,
\epsilon_1^{}\!-\epsilon_{N-1}^{},\underline{\epsilon_1^{}\!-
\epsilon_N^{}},\epsilon_{N-1}^{}\!-\epsilon_{N}^{},\ldots,\epsilon_2^{}\!-
\epsilon_{N}^{})~,&
\\[5pt]
&(\epsilon_2^{}\!-\epsilon_3^{},\epsilon_2^{}\!-\epsilon_4^{},\ldots,
\epsilon_2^{}\!-\epsilon_{N-2}^{},\underline{\epsilon_2^{}\!-
\epsilon_{N-1}^{}},\epsilon_{N-2}^{}\!-\epsilon_{N-1}^{},\ldots,
\epsilon_3^{}\!-\epsilon_{N-1}^{})~,&
\\[1pt]
&\ldots\ldots\ldots\ldots\ldots\ldots\ldots\ldots
\ldots\ldots\ldots\ldots\ldots\ldots\ldots\ldots\ldots\ldots,&
\end{array}
\label{jr23}
\end{equation}
The underlined roots are maximal on their lines. Each line of
(\ref{jr23}) is corresponding to an extended Jordanian classical
$r$-matrix of maximal order with its parameter deformation. Thus
we have the set of such $r$-matrices:
\begin{equation}
\begin{array}{rcl}
r_{1}^{}\left(\xi_{1}^{}\right)\!\!&:=\!\!&\xi_{1}^{}\biggl(\displaystyle{
\frac{1}{2}}(e_{1-1}^{}-e_{N-N}^{})\wedge e_{1-N}^{}+\sum_{i=2}^{N-1}
(-1)^{\deg e_{1-i}\deg e_{i-N}}e_{1-i}\wedge e_{i-N}\biggr),
\\[7pt]
r_{2}^{}\left(\xi_{2}^{}\right)\!\!&:=\!\!&\xi_{2}^{}\biggl(\displaystyle{
\frac{1}{2}}(e_{2-2}^{}\!-\!e_{N\!-\!1,-N\!+\!1}^{})\wedge
e_{2,-N\!+\!1}^{}+\sum_{i=2}^{N-2}
(-1)^{\deg e_{2-i}\deg e_{i,-N\!+\!1}}e_{2-i}\wedge e_{i,-N\!+\!1}\biggr),
\\[1pt]
\ldots\!\!&\!\!&\ldots\ldots\ldots\ldots\ldots\ldots\ldots\ldots
\ldots\ldots\ldots\ldots\ldots\ldots\ldots\ldots\ldots\ldots\ldots\ldots~,
\label{jr24}
\end{array}
\end{equation}
where the root vectors $e_{i-k}:=e_{\epsilon_i-\epsilon_k}$
($i<k$) are chosen such that they satisfy the relations
(\ref{jr12}). The resulting maximal $r$-matrix is the sum of these
matrices:
\begin{equation}
r_{1,2,\ldots,[N/2]}^{}
\left(\xi_{1},\xi_{2}^{},\cdots,\xi_{[N/2]}^{}\right)=r_{1}^{}(\xi_{1}^{})+
r_{2}^{}(\xi_{2}^{})+\cdots+r_{[N/2]}^{}(\xi_{[N/2]}^{})~.
\label{jr25}
\end{equation}
\begin{proposition} The elements of the subalgebra
$\mathfrak{gl}(m-i|n-i)$ co-commute with the $r$-matrix
\begin{equation}
r_{1,2,\ldots,i}^{}\left(\xi_{1}^{},\xi_{2}^{},\ldots,\xi_{i}^{}\right):=
r_{1}^{}\left(\xi_{1}^{}\right)+r_{2}^{}\left(\xi_{2}^{}\right)+\cdots+
r_{i}^{}\left(\xi_{i}^{}\right)~,
\label{jr26}
\end{equation}
i.e.
\begin{equation}
[x\otimes1+1\otimes x,\,r_{1,2,\ldots,i}^{}
\left(\xi_{1}^{},\xi_{2}^{},\ldots,\xi_{i}^{}\right)]=0\quad
(\forall x\in\mathfrak{gl}(m-i|n-i)).
\label{jr27}
\end{equation}
\end{proposition}
Thus the constructed extended Jordanian $r$-matrices
$r_{1}^{}\left(\xi_{1}^{}\right)$, $r_{2}^{}\left(\xi_{2}^{}\right)$,
$\ldots$, $r_{[N/2]}^{}(\xi_{[N/2]}^{})$
are associated with the following reduction chain
\begin{equation}
\mathfrak{gl}(m|n)\supset\mathfrak{gl}(m-1|n-1)\supset
\mathfrak{gl}(m-2|n-2)^{}\cdots\supset\mathfrak{gl}(k)
\quad(k=3\;{\rm or}\;2|1,{\rm or}\;2,{\rm or}\; 1|1)~.
\label{jr28}
\end{equation}
{{\bf Example 4. Maximal classical r-matrix of Jordanian type for
the Lie superalgebra $C(n)\simeq\mathfrak{osp}(1|2n)$}. In the
terms of the orthonormalized basis $\epsilon_i^{}$
($i=1,2,\ldots,n$) the systems of positive roots, $\Delta_{+}^{}$,
for $\mathfrak{osp}(1|2n)$ are given as follows:
\begin{eqnarray}
\Delta_{+}(\mathfrak{osp}(1|2n))\!\!\!\!&=\!\!&
\{\epsilon_i^{}\pm\epsilon_j^{},\,
\epsilon_k^{},\,2\epsilon_k^{}\,|\;1\le i<j\le n;\,
k=1,2,\ldots,n\}~. \label{jr29}
\end{eqnarray}
Let us write down the positive root system
$\Delta_{+}^{}(\mathfrak{osp}(1|2n)$
in the following ordering
\begin{equation}
\begin{array}{rcl}
&(\epsilon_1^{}\!-\epsilon_2^{},\epsilon_1^{}\!-\epsilon_3^{},\ldots,
\epsilon_1^{}\!-\epsilon_n^{},\epsilon_1^{},\underline{2\epsilon_1^{}},
\epsilon_1^{}\!+\epsilon_n^{},\ldots,\epsilon_{1}^{}\!+\epsilon_{3}^{},
\epsilon_{1}^{}\!+\epsilon_{2}^{}),&
\\[5pt]
&(\epsilon_2^{}\!-\epsilon_3^{},\epsilon_2^{}\!-\epsilon_4^{}\ldots,
\epsilon_2^{}\!-\epsilon_n^{},\epsilon_2^{},
\underline{2\epsilon_2^{}},\epsilon_2^{}\!+\epsilon_n^{},\ldots,
\epsilon_{n-1}^{}\!+\epsilon_{4}^{},\epsilon_{2}^{}\!+\epsilon_{3}^{}),&
\\[1pt]
&\ldots\ldots\ldots\ldots\ldots\ldots\ldots\ldots\ldots\ldots\ldots
\ldots\ldots\ldots\ldots\ldots\ldots&
\\[1pt]
&\epsilon_{n}^{},\underline{2\epsilon_{n}^{}}~.&
\end{array}
\label{jr30}
\end{equation}
Each line of (\ref{jr30}) is corresponding to a extended Jordanian
classical $r$-matrix of maximal order with its parameter
deformation, and we have the set of such $r$-matrices:
\begin{equation}
\begin{array}{rcl}
r_{1}^{}\left(\xi_{1}^{}\right)\!\!&:=\!\!&\xi_{1}^{}\biggl
(\displaystyle{\frac{1}{2}}e_{1-1}^{}\wedge
e_{11}^{}-e_{1}^{}\wedge e_{1}^{}+\sum_{i=2}^{n}e_{1-i}\wedge
e_{1i}\biggr)~,
\\
r_{2}^{}\left(\xi_{2}^{}\right)\!\!&:=\!\!&\xi_{2}^{}\biggl
(\displaystyle{\frac{1}{2}}e_{2-2}\wedge e_{22}^{}-e_{2}^{}\wedge
e_{2}^{}+\sum_{i=3}^{n}e_{2-i}\wedge e_{2i}\biggr)~,
\\[1pt]
\ldots&&\ldots\ldots\ldots\ldots\ldots\ldots\ldots\ldots\ldots
\ldots\ldots\ldots
\\[1pt]
r_{n}^{}\left(\xi_{n}^{}\right)\!\!&:=\!\!&\xi_{n}^{}\,\biggl(
\displaystyle{\frac{1}{2}}e_{n-n}\wedge e_{nn}^{}-e_{n}^{}\wedge
e_{n}^{}\biggr)~,
\end{array}
\label{jr31}
\end{equation}
where the root vectors $e_{i-k}^{}:=e_{\epsilon_i-\epsilon_k}$,
$e_{ik}^{}:=e_{\epsilon_i+\epsilon_k}$,
$e_{i}^{}:=e_{\epsilon_i}^{}$ are chosen such that they satisfy
the relations (\ref{jr12}). The resulting maximal $r$-matrix is
the sum of these matrices:
\begin{equation}
r_{1,2,\ldots,n}^{}\left(\xi_{1},\xi_{2}^{},\cdots,\xi_{n}^{}\right)=
r_{1}^{}(\xi_{1}^{})+r_{2}^{}(\xi_{2}^{})+\cdots+r_{n}^{}(\xi_{n}^{})~.
\label{jr32}
\end{equation}
\begin{proposition} The elements of the subalgebra $\mathfrak{osp}(1|2n-2i)$
co-commute with the $r$-matrix
\begin{equation}
r_{1,2,\ldots,i}^{}\left(\xi_{1}^{},\xi_{2}^{},\ldots,\xi_{i}^{}\right):=
r_{1}^{}\left(\xi_{1}^{}\right)+r_{2}^{}\left(\xi_{2}^{}\right)+\cdots+
r_{i}^{}\left(\xi_{i}^{}\right)~,
\label{jr33}
\end{equation}
i.e.
\begin{equation}
[x\otimes1+1\otimes x,\,r_{1,2,\ldots,i}^{}
\left(\xi_{1}^{},\xi_{2}^{},\ldots,\xi_{i}^{}\right)]=0\qquad
\bigl(\forall x\in\mathfrak{osp}(1|2n-2i)\bigr)~. \label{jr34}
\end{equation}
\end{proposition}
Thus the constructed extended Jordanian classical $r$-matrices
$r_{1}^{}\left(\xi_{1}^{}\right)$, $r_{2}^{}\left(\xi_{2}^{}\right)$,$\ldots$,
$r_{n}^{}(\xi_{n}^{})$ are associated with the following reduction chain
\begin{equation}
\mathfrak{osp}(1|2n)\supset\mathfrak{osp}(1|2n-2)\supset
\mathfrak{osp}(1|2n-4)^{}\cdots\supset\mathfrak{osp}(1|2)~.
\label{jr35}
\end{equation}
{\bf Example 5. Maximal classical r-matrix of Jordanian type for the
Lie superalgebras $B(m|n)\simeq\mathfrak{osp}(2m+1|2n)$ and
$D(m|n)\simeq\mathfrak{osp}(2m|2n)$}. In the terms of the orthonormalized
basis $\epsilon_i^{}$ ($i=1,2,\ldots,m+n$) the systems of positive roots
for $\mathfrak{osp}(2m+1|2n)$ are given as follows:
\begin{equation}
\Delta_{+}^{}(\mathfrak{osp}(2m+1|2n))=
\{\epsilon_i^{}\pm\epsilon_j^{},\epsilon_k^{},\epsilon_{2l}^{}\,
|\,1\le i<j\le N;\,1\le k\le N;\,m+1\le l\le N\}. \label{jr36}
\end{equation}
where $N:=m+n$.
Let us write down the positive root system in the following ordering
\begin{equation}
\begin{array}{rcl}
&(\epsilon_1^{}\!-\epsilon_2^{}),(\epsilon_1^{}\!-\epsilon_3^{},
\ldots,\epsilon_1^{}\!-\epsilon_N^{},\epsilon_1^{},\epsilon_1^{}\!+
\epsilon_N^{},\ldots,\epsilon_{1}^{}\!+\epsilon_{3}^{},
\underline{\epsilon_{1}^{}\!+\epsilon_{2}^{}},\epsilon_2^{}\!-
\epsilon_3^{},
\phantom{aaaaaaaaa}&
\\[2pt]
&\phantom{aaaaaaaaaaaaaaaaaa}\epsilon_2^{}\!-\epsilon_4^{},\ldots,
\epsilon_2^{}\!-\epsilon_{N}^{},\epsilon_2^{},\epsilon_2^{}\!+
\epsilon_{N}^{},\ldots,\epsilon_{2}^{}\!+\epsilon_{4}^{},
\epsilon_{2}^{}\!+\epsilon_{3}^{}),&
\\[7pt]
&(\epsilon_3^{}\!-\epsilon_4^{}),(\epsilon_3^{}\!-\epsilon_5^{},
\ldots,\epsilon_3^{}\!-\epsilon_{N}^{},\epsilon_3^{},\epsilon_3^{}\!+
\epsilon_{N}^{}, \ldots,\epsilon_{3}^{}\!+\epsilon_{5}^{},
\underline{\epsilon_{3}^{}\!+\epsilon_{4}^{}},
\epsilon_4^{}\!-\epsilon_5^{}, \phantom{aaaaaaaaa}&
\\[2pt]
&\phantom{aaaaaaaaaaaaaaaaaa}\epsilon_4^{}\!-\epsilon_6^{},\ldots,
\epsilon_4^{}\!-\epsilon_{N}^{},\epsilon_3^{},\epsilon_4^{}\!+
\epsilon_{N}^{},\ldots,\epsilon_{2}^{}\!+\epsilon_{4}^{},
\epsilon_{2}^{}\!+\epsilon_{3}^{}),&
\\[1pt]
&\ldots\ldots\ldots\ldots\ldots\ldots\ldots\ldots\ldots\ldots
\ldots\ldots\ldots\ldots\ldots\ldots\ldots\ldots\ldots\ldots&
\\[1pt]
&(\epsilon_{m'}^{}\!-\epsilon_{m'+1}^{},\ldots,\epsilon_{m'}^{}\!-
\epsilon_N^{},\epsilon_{m'}^{},\underline{2\epsilon_{m'}^{}},
\epsilon_{m'}^{}\!+\epsilon_N^{},\ldots,\epsilon_{m}^{}\!+
\epsilon_{m'+1}^{}),&
\\[1pt]
&\ldots\ldots\ldots\ldots\ldots\ldots\ldots\ldots\ldots
\ldots\ldots\ldots\ldots\ldots\ldots \ldots\ldots\ldots&
\\[1pt]
&\epsilon_{N}^{},\underline{2\epsilon_{N}^{}}~,&
\end{array}
\label{jr37}
\end{equation}
where $m'=m+1$ if $m$ is an even positive integer and $m'=m+2$ if
$m$ is an odd positive integer. Each set of roots in the brackets
$(\ldots)$ is corresponding to an extended Jordanian classical
$r$-matrix of maximal order with its parameter deformation:
\begin{equation}
\begin{array}{rcl}
r_{1}^{}\left(\xi_{1}^{},\xi_{1}'\right)\!\!&:=\!\!&\xi_{1}^{}
\biggl(\displaystyle{\frac{1}{2}}(e_{1-1}^{}+e_{2-2}^{})\wedge
e_{12}^{}+ (-1)^{\deg e_{1}^{}\deg e_{2}^{}}e_{1}^{}\wedge
e_{2}^{}
\\[9pt]
&&+\sum\limits_{i=3}^{N}(-1)^{\deg e_{1-i}^{}\deg
e_{2i}^{}}e_{1-i} \wedge e_{2i}+\sum\limits_{i=3}^{N}(-1)^{\deg
e_{1i}^{}\deg e_{2-i}^{}}e_{1i} \wedge e_{2-i}\biggr),
\\[9pt]
&&+\xi_{1}'\displaystyle{\frac{1}{2}}(e_{1-1}^{}-e_{2-2}^{})
\wedge e_{1-2}^{}~,
\\[15pt]
r_{2}^{}\left(\xi_{2}^{},\xi_{2}'\right)\!\!&:=\!\!&\xi_{2}^{}
\biggl(\displaystyle{\frac{1}{2}}(e_{3-3}^{}+e_{4-4}^{})\wedge
e_{34}^{}+ (-1)^{\deg e_{3}^{}\deg e_{4}^{}}e_{3}^{}\wedge
e_{4}^{}+
\\[9pt]
&&+\sum\limits_{i=5}^{N}(-1)^{\deg e_{3-i}^{}\deg e_{4i}^{}}e_{3-i}\wedge
e_{4i}+\sum\limits_{i=5}^{N}(-1)^{\deg e_{3i}^{}\deg e_{4-i}^{}}e_{3i}
\wedge e_{4-i}\biggr)
\\[9pt]
&&+\xi_{2}'\displaystyle{\frac{1}{2}}(e_{3-3}^{}-e_{4-4}^{})
\wedge e_{3-4}^{}~,
\\[1pt]
\ldots&&\ldots\ldots\ldots\ldots\ldots\ldots\ldots\ldots
\ldots\ldots\ldots\ldots\ldots\ldots\ldots\ldots\ldots\ldots
\\[1pt]
r_{m'}^{}\left(\xi_{m'}^{}\right)\!\!&:=\!\!&\xi_{m'}^{}
\biggl(\displaystyle{\frac{1}{2}}e_{m'-m'}\wedge e_{m'm'}^{}-
e_{m'}^{}\wedge e_{m'}^{}\!+\!\!\sum_{i=m'+1}^{N}e_{m'-i}
\wedge e_{m'i}\biggr),
\\[1pt]
\ldots&&\ldots\ldots\ldots\ldots\ldots\ldots\ldots\ldots
\ldots\ldots\ldots\ldots\ldots\ldots\ldots\ldots\ldots\ldots
\\[1pt]
r_{N}^{}\left(\xi_{N}^{}\right)\!\!&:=\!\!&\xi_{N}^{}\,\biggl(
\displaystyle{\frac{1}{2}}e_{N-N}\wedge e_{NN}^{}-e_{N}^{}\wedge
e_{N}^{}\biggr)~.
\end{array}
\label{jr38}
\end{equation}
The resulting maximal $r$-matrix is the sum of these matrices:
\begin{equation}
\begin{array}{lcr}
&&r_{1,2,\ldots,N}^{}\left(\xi_{1},\xi_{1}',\xi_{2}^{},\xi_{2}',
\ldots,\xi_{m'},\ldots,\xi_{N}^{}\right)\;=\;
\phantom{aaaaaaaaaaaaaaaaaaaaaaaaaaaaaaa}
\\[5pt]
&&=\;r_{1}^{}(\xi_{1}^{},\xi_{1}')+r_{2}^{}(\xi_{2}^{},\xi_{2}')+
\cdots+r_{m'}^{}(\xi_{m'}^{})+\cdots+r_{N}^{}(\xi_{N}^{})~.
\end{array}
\label{jr39}
\end{equation}
\begin{proposition} The elements of the subalgebra
$\mathfrak{so}(3)\oplus\mathfrak{osp}(2(m-i)-1|2n)$
co-commute with the $r$-matrix
\begin{equation}
r_{1,2,\ldots,i}^{}\left(\xi_{1}^{},\xi_{1}',\xi_{2}^{},\xi_{2}',
\ldots,\xi_{i}^{},\xi_{i}',\right):=r_{1}^{}\left(\xi_{1}^{},
\xi_{1}'\right)+r_{2}^{}\left(\xi_{2}^{},\xi_{2}'\right)+\cdots+
r_{i}^{}\left(\xi_{i}^{},\xi_{i}'\right)~, \label{jr40}
\end{equation}
i.e.
\begin{equation}
[x\otimes1+1\otimes x,\,r_{1,2,\ldots,i}^{}\left(\xi_{1}^{},
\xi_{2}^{},\xi_{2}',\ldots,\xi_{i}^{},\xi_{i}'\right)]=0\quad
(\forall x\in\mathfrak{so}(3)\oplus\mathfrak{osp}(2(m-i)-1|2n))~.
\label{jr41}
\end{equation}
\end{proposition}
Thus the constructed extended Jordanian classical $r$-matrices
$r_{1}^{}\left(\xi_{1}^{},\xi_{1}'\right)$,
$r_{2}^{}\left(\xi_{2}^{},\xi_{2}'\right),\ldots,
,r_{N}^{}(\xi_{N}^{})$ are associated with the following reduction
chain
\begin{equation}
\mathfrak{osp}(2m+1|2n)\supset\mathfrak{so}(3)\oplus\mathfrak{osp}(2m-3|2n)
\supset\ldots\supset\mathfrak{osp}(1|2n')\supset\cdots
\supset\mathfrak{osp}(1|2) \label{jr42}
\end{equation}
where $n'=n$ if $m$ is an even positive integer, and $n'=n-1$ if
$m$ is an odd positive integer.

We obtain the results for the Lie superalgebra
$D(m|n)\simeq\mathfrak{osp}(2m|2n)$} if we remove all roots
$\varepsilon_i$ and the root vectors $e_i$ $(i=1, \ldots, N)$
in the formulas (\ref{jr37}), (\ref{jr38}).

\setcounter{equation}{0}
\section{Chains of extended Jordanian twists}

Basic elements of Drinfeld`s theory of twisting quantization for Hopf
algebras \cite{D} is easy generalized to the case of Hopf superalgebras.
Indeed, formulae describing twist quantization for Lie algebras and Lie
superalgebras are the same, provided in the second case the grading
is properly taken into account.

Let $\mathfrak{U}:=\mathfrak{U}(m,\Delta,S,\epsilon)$ be a Hopf
superalgebra with graded structure
($\mathfrak{U}=\mathfrak{U}_0\oplus\mathfrak{U}_1$) and with the
grading preserving operations, a multiplication $m:\mathfrak{U}
\otimes\mathfrak{U}\rightarrow\mathfrak{U}$, a coproduct $\Delta:
\mathfrak{U}\rightarrow\mathfrak{U}\otimes\mathfrak{U}$, an
antipode $S:\mathfrak{U}\to\mathfrak{U}$, and a counit
$\epsilon:\mathfrak{U} \rightarrow\mathbb C$. Let there exists an
invertible even element $F=\sum_i f^{(1)}_i\otimes f^{(2)}_i$ of
some extension of $\mathfrak{U}\otimes\mathfrak{U}$, such that it
satisfies the cocycle equation
\begin{equation}
F^{12}(\Delta\otimes{\rm id})(F)=F^{23}({\rm id}\otimes\Delta)(F)~,
\label{jt1}
\end{equation}
and the "unital" normalization condition
\begin{equation}
(\epsilon \otimes{\rm id})(F)=({\rm id}\otimes\epsilon )(F)=1~.
\label{jt2}
\end{equation}
Then the new Hopf superalgebra $\mathfrak{U}^{(F)}:=
\mathfrak{U}^{(F)}(m,\Delta^{(F)},S^{(F)},\epsilon$) with the same
multiplication $m$ and the counit mapping $\epsilon$ but with the
twisted coproduct and antipode
\begin{equation}
\Delta^{(F)}(a)=F\Delta(a)F^{-1},\quad\;S^{(F)}=u\,S(a)u^{-1},
\quad u=\sum_i f^{(1)}_{i}S(f^{(2)}_i)\;\quad(\forall a\in\mathfrak{A})
\label{jt3}
\end{equation}
is called the twisted Hopf superalgebra or twist quantization of
the Hopf superalgebra $\mathfrak{U}$.

The Hopf superalgebra $\mathfrak{U}$ is called quasitriangular if it
has an additional invertible element (universal $R$-matrix) $R$  which
relates the coproduct $\Delta$ with its opposite coproduct $\tilde{\Delta}$
by the similarity transformation
\begin{equation}
\tilde{\Delta}(a)=R\,\Delta(a)R^{-1} \qquad(\forall a\in
\mathfrak{U})~, \label{jt4}
\end{equation}
and $R$ satisfies the quasitriangularity conditions
\begin{equation}
(\Delta\otimes {\rm id})(R)=R^{13}R^{23}~,
\qquad ({\rm id}\otimes\Delta)(R)=R^{13}R^{12}~.
\label{jt5}
\end{equation}
The twisted ("quantized") Hopf algebra $\mathfrak{U}^{(F)}$ is also
quasitriangular with the universal $R$-matrix $R^{(F)}$ defined as
follows
\begin{equation}
R^{(F)} = F^{21}R\,F^{-1}~,
\label{jt6}
\end{equation}
where $F^{21}=\sum_i f^{(2)}_i\!\otimes f^{(1)}_{i}$ provided
$F=\sum_i f^{(1)}_i\!\otimes f^{(2)}_i$. For the non-deformed,
classical case $\mathfrak{U}=U(\mathfrak{g})$, where $\mathfrak{g}$
is a simple Lie superalgebra, the universal $R$-matrix is trivial, $R=1$.

The twisting two-tensor $F_{\theta,N}^{}(\xi)$ corresponding to
the $r$-matrix (\ref{jr14}) has the form
\begin{equation}
F_{\theta,N}^{}(\xi)\;=\;\mathfrak{F}_{N}^{}(\xi)F_{\!J}^{}
\label{jt7}
\end{equation}
where the two-tensor $F_{\!J}^{}$ is the Jordanian twist
corresponding to the Jordanian $r$-matrix (\ref{jr2}) and
$\mathfrak{F}_{N}^{}$ is extension of the Jordanian twist. These
two-tensors are given by the formulas
\begin{equation}
F_{\!J}^{}\;= \;\exp(2h_{\theta}^{}\otimes\sigma_{\theta}^{})~,
\label{jt8}
\end{equation}
\begin{equation}
\begin{array}{rcl}
\mathfrak{F}_{N}^{}(\xi)\!\!&=\!\!&\bigg(\prod\limits_{i=1}^{N'}
\exp\Big(\xi(-1)^{\deg e_{\gamma_i}^{}\deg e_{\gamma_{-i}}^{}}
e_{\gamma_i}^{}\otimes e_{\gamma_{-i}}^{}
\exp(-2t_{\gamma_i}\sigma_{\theta}^{})\Big)\bigg)\mathfrak{F}_S^{}
\\
&=\!\!&\exp\Big(\xi\sum\limits_{i=1}^{N'} (-1)^{\deg
e_{\gamma_i}^{}\deg e_{\gamma_{-i}}^{}}e_{\gamma_i}^{}\otimes
e_{\gamma_{-i}}^{}\exp(-2t_{\gamma_i}\sigma_{\theta}^{})\Big)
\,\mathfrak{F}_S^{}~, \label{jt9}
\end{array}
\end{equation}
where
\begin{equation}
\mathfrak{F}_S^{}=\Biggl(1-\xi\frac{e_{\theta/2}^{}}{\exp\sigma_{\theta}+1}
\otimes\frac{e_{\theta/2}^{}}{\exp\sigma_{\theta}+1}\Biggr)
\sqrt{\frac{(e^\sigma+1)\otimes(e^\sigma+1)} {2(e^\sigma\!\otimes
e^\sigma+1)}}~, \label{jt10}
\end{equation}
if $\theta/2$ is a root (i.e. $\theta/2=\gamma_i=\gamma_{-i}$ for
some $i$), $e_{\theta/2}^2=e_{\theta}^{}$, $N'=N-1$, and
\begin{equation}
\mathfrak{F}_S^{}=1~, \label{jt11}
\end{equation}
if $\theta/2$ is not any root, $N'=N$. Moreover
\begin{equation}
\deg(\xi)=\deg(e_{\theta}^{})=\deg(e_{\gamma_i})+\deg(e_{\gamma_{-i}})
\quad(\mathop{\rm mod}2)~, \label{jt12}
\end{equation}
\begin{equation}
\sigma_{\theta}^{}:=\frac{1}{2}\ln(1+\xi e_{\theta}^{})~.
\label{jt13}
\end{equation}
It should be noted that if the root vector $e_{\theta}^{}$ is odd
then $\sigma_{\theta}^{}=\frac{1}{2}\xi e_{\theta}^{}$.

The twisted coproduct $\Delta_{\xi}(\,\cdot\,):=
F_{\theta,N}^{}(\xi)\Delta(\,\cdot\,)F_{\theta,N}^{-1}(\xi)$ and
the corresponding antipode $S_{\xi}^{}$ for elements in
(\ref{jr12}) are given by the formulas
\begin{eqnarray}
\Delta_{\xi}(\exp(\pm\sigma_{\theta}^{}))\!\!&=\!\!&
\exp(\pm\sigma_{\theta}^{})\otimes\exp(\pm\sigma_{\theta}^{}),\quad
\Delta_{\xi}(e_{\theta/2}^{})\;=\;e_{\theta/2}^{}\otimes1+
\exp(\sigma_{\theta}^{})\otimes e_{\theta/2}^{}~,\quad
\label{jt14}
\\[5pt]
\Delta_{\xi}(h_{\theta}^{})\!\!&=\!\!&
h_{\theta}^{}\otimes\exp(-2\sigma_{\theta}^{})+1\otimes h_{\theta}^{}+
\frac{\xi}{4}e_{\theta/2}^{}\exp(-\sigma_{\theta}^{})\otimes
e_{\theta/2}^{}\exp(-2\sigma_{\theta}^{})\nonumber
\\[-4pt]
&&-\xi\sum\limits_{i=1}^{N'}(-1)^{\deg e_{\gamma_i}^{}\deg
e_{\gamma_{-i}}^{}} e_{\gamma_i}^{}\otimes e_{\gamma_{-i}}^{}
\exp(-2(t_{\gamma_i}^{}+1)\sigma_{\theta}^{})~, \label{jt15}
\\[5pt]
\Delta_{\xi}(e_{\gamma_i}^{})\!\!&=\!\!&
e_{\gamma_i}^{}\otimes\exp(-2t_{\gamma_i}\sigma_{\theta}^{})+
1\otimes e_{\gamma_i}^{}~, \label{jt16}
\\[5pt]
\Delta_{\xi}(e_{\gamma_{-i}}^{})\!\!&=\!\!&
e_{\gamma_{-i}}^{}\otimes\exp(2t_{\gamma_i}\sigma_{\theta}^{})+
\exp(2\sigma_{\theta}^{})\otimes e_{\gamma_{-i}}^{}~, \label{jt17}
\end{eqnarray}
\begin{eqnarray}
S_{\xi}(\exp(\pm\sigma_{\theta}^{}))\!\!&=\!\!&
\exp(\mp\sigma_{\theta}^{})~,\qquad\quad
S_{\xi}(e_{\theta/2}^{})\;=
\;-e_{\theta/2}^{}\exp(-\sigma_{\theta}^{})~, \label{jt18}
\\[1pt]
S_{\xi}(h_{\theta}^{})\!\!&=\!\!&-h_{\theta}^{}\exp(2\sigma_{\theta}^{})+
\frac{1}{4}\Bigl(\exp(2\sigma_{\theta}^{})-1\Bigr)-\xi\sum_{i=1}^{N'}
(-1)^{\deg e_{\gamma_i}^{}\deg e_{\gamma_{-i}}^{}}
e_{\gamma_i}^{}e_{\gamma_{-i}}^{}, \label{jt19}
\\[1pt]
S_{\xi}(e_{\gamma_i}^{})\!\!&=\!\!&-e_{\gamma_i}^{}
\exp(2t_{\gamma_i}^{}\sigma_{\theta}^{})~,\quad
S_{\xi}(e_{\gamma_{-i}}^{})\;=\;-e_{\gamma_{-i}}^{}
\exp(-2(t_{\gamma_i}^{}\!+1)\sigma_{\theta}^{})~. \label{jt20}
\end{eqnarray}
If $\theta/2$ is not any root, the third term in (\ref{jt15}) and
the second term in (\ref{jt19}) should be removed.
\begin{proposition}
(i) Let $\mathfrak{g}'$ be a subalgebra of $\mathfrak{g}$, which
co-commutes with the $r$-matrix (\ref{jr14}),
$\delta_{\xi}(\mathfrak{g}')=0$ (see (\ref{jr18})). If an
invertible element $w_{\xi}^{}$ of some extension of
$U(\mathfrak{g})$ satisfies the equations
\begin{equation}
[\Delta(x),(w_{\xi}^{-1}\otimes
w_{\xi}^{-1})\Delta_{\xi}(w_{\xi}^{})F_{\theta,N}^{}(\xi)]\;=
\;0\quad (\forall x\in\mathfrak{g}'),\quad w_{\xi}^{}\equiv1\;
\mathop{\rm mod} \xi,\quad\varepsilon(w_{\xi}^{})=1,\label{jt21}
\end{equation}
then the automorphism $w_{\xi}^{}x w_{\xi}^{-1}$ simplifies (makes
trivial) the twisted coproduct $\Delta_{\xi}(\,\cdot\,)$ in
$\mathfrak{g}'$:
\begin{equation}
\Delta_{\xi}(w_{\xi}^{}x w_{\xi}^{-1})\;:=\;w_{\xi}^{}x
w_{\xi}^{-1}\otimes1+1\otimes w_{\xi}^{}x w_{\xi}^{-1}~,\quad
x\in\mathfrak{g}'.\label{jt22}
\end{equation}
(ii) The element
$w_{\xi}^{}\equiv\sqrt{u(\mathfrak{F}_{N}^{}(\xi))}$ satisfies the
equations (\ref{jt21}), where $u(\mathfrak{F}_{N}^{}(\xi))$ is the
Hopf "folding" of the two-tensor (\ref{jt9}):
\begin{eqnarray}
u(\mathfrak{F}_{N}^{}(\xi))\!\!&=\!\!&((S_{\xi}^{}\otimes
\mathop{\rm Id})\mathfrak{F}_{N}^{}(\xi))\circ1\,= \,
((\mathop{\rm Id}\otimes S_{J}^{} )\mathfrak{F}_{N}^{}(\xi)) \circ
1 \label{jt23}
\\[5pt]
&=&\Biggl(\prod\limits_{i=1}^{N'}\biggl(\sum_{n=0}^{\infty}
\frac{(-\xi)^n}{n!}(-1)^{n\deg e_{\gamma_i}^{}\deg
e_{\gamma_{-i}}^{}} e_{\gamma_i}^{n}
e_{\gamma_{-i}}^{n})\biggr)\Biggr)u_S^{}\label{jt24}
\\[5pt]
&=&\exp\bigg(\frac{-2\xi\sigma_{\theta}}{\exp(2\sigma_{\theta})-1}
\sum_{i=1}^{N'}(-1)^{\deg e_{\gamma_i}^{} \deg
e_{\gamma_i}^{}}e_{\gamma_i}^{}e_{\gamma_{-i}}^{}\Big)u_S~.
\label{jt25}
\end{eqnarray}
Here $u_S^{}$ is the folding of the super-tensor
$\mathfrak{F}_S^{}$: $u_S=\exp(\frac{1}{2}\sigma)$ if $\theta/2$
is a root, and $u_S=1$ if $\theta/2$ is not any root; $S_{J}^{}$
is the antipode after the Jordanian twist (\ref{jt8}); the
operation $"\circ"$ means $(a\otimes b)\circ x=axb$. The inverse
element $u^{-1}(\mathfrak{F}_{N}^{}(\xi))$ is given the following
explicit formula
\begin{eqnarray}
u^{-1}(\mathfrak{F}_{N}^{}(\xi))\!\!&=\!\!&((\mathop{\rm Id}
\otimes S_{\xi}^{} )\mathfrak{F}_{N}^{-1}(\xi))\circ 1\,=\,
((S_{J}^{}\otimes \mathop{\rm Id})\mathfrak{F}_{N}^{-1}(\xi))
\circ 1 \label{jt26}
\\[5pt]
&=&\Biggl(\prod\limits_{i=1}^{N'}\biggl(\sum_{n=0}^{\infty}
\frac{(-\xi)^n\exp(-2n\sigma)}{n!}(-1)^{n\deg e_{\gamma_i}^{}\deg
e_{\gamma_{-i}}^{}}e_{\gamma_i}^{n}e_{\gamma_{-i}}^{n})\biggr)\Biggr)
u_S^{-1}\phantom{aa}\label{jt27}
\\[5pt]
&=&\exp\bigg(\frac{2\xi\sigma_{\theta}}{\exp(2\sigma_{\theta})-1}
\sum_{i=1}^{N'}(-1)^{\deg e_{\gamma_i}^{}\deg e_{\gamma_i}^{}}
e_{\gamma_i}^{}e_{\gamma_{-i}}^{}\Big)u_S^{-1}~.\label{jt28}
\end{eqnarray}
Moreover
\begin{eqnarray}
\sqrt{u(\mathfrak{F}_{N}^{}(\xi))}\!\!&=\!\!&
\Biggl(\prod\limits_{i=1}^{N'}\biggl(\sum_{n=0}^{\infty}
\frac{(-\xi)^{n}(-1)^{n\deg e_{\gamma_i}^{}\deg
e_{\gamma_{-i}}^{}}}{n!(\exp(\sigma_{\theta}^{})+1)^n}\;
e_{\gamma_i}^{n} e_{\gamma_{-i}}^{n}\Biggr)
\sqrt{u_S^{}}\label{jt29}
\\[5pt]
&=&\exp\bigg(\frac{-\xi\sigma_{\theta}}{\exp(2\sigma_{\theta})-1}
\sum_{i=1}^{N'}(-1)^{\deg e_{\gamma_i}^{} \deg
e_{\gamma_i}^{}}e_{\gamma_i}^{}e_{\gamma_{-i}}^{}\Big)
\sqrt{u_S}~, \label{jt30}
\\[7pt]
\sqrt{u^{-1}(\mathfrak{F}_{N}^{}(\xi))}\!\!&=\!\!&
\Biggl(\prod\limits_{i=1}^{N'}\biggl(\sum_{n=0}^{\infty}
\frac{\xi^{n}\exp(-n\sigma)}{n!(\exp(\sigma_{\theta}^{})+1)^n}\,
(-1)^{n\deg e_{\gamma_i}^{}\deg
e_{\gamma_{-i}}^{}}e_{\gamma_i}^{n}
e_{\gamma_{-i}}^{n}\Biggr)\sqrt{u_S^{-1}}\label{jt31}
\\[5pt]
&=&\exp\bigg(\frac{\xi\sigma_{\theta}}{\exp(2\sigma_{\theta})-1}
\sum_{i=1}^{N'}(-1)^{\deg e_{\gamma_i}^{} \deg
e_{\gamma_i}^{}}e_{\gamma_i}^{}e_{\gamma_{-i}}^{}\Big)
\sqrt{u_S^{-1}}~. \label{jt32}
\end{eqnarray}
\end{proposition}
{\it Remarks.} The formula (\ref{jt30}) for the case
$\mathfrak{g}=\mathfrak{sp}(2n)$ was found in \cite{AKL} by
fitting.

With the help of the elements
$w_{\xi}^{}=\sqrt{u(\mathfrak{F}_{N}^{}(\xi))}$,
$w_{\xi_1^{}}^{}=\sqrt{u(\mathfrak{F}_{N_1}^{}(\xi_1^{}))},\;\ldots\,,
w_{\xi_{k}}^{}=\sqrt{u(\mathfrak{F}_{N_k}^{}(\xi_{k}^{}))}$ the
total twist chain corresponding to the $r$-matrix (\ref{jr21}) can
be presented as follows
\begin{equation}
\begin{array}{rcl}
F_{\theta,N;\theta_1^{},N_1^{};\ldots;\theta_k^{},N_k^{}}^{}
(\xi,\xi_1^{},\ldots,\xi_k^{})\!\!&=\!\!&
F_{\theta_k^{},N_k^{}}^{}(\xi,\xi_1^{},
\ldots,\xi_{k-1}^{};\xi_k^{})\cdots
\\[9pt]
&&\times F_{\theta_{2}^{},N_{2}^{}}^{}(\xi,\xi_1^{};\xi_2)
F_{\theta_{1}^{},N_{1}^{}}^{}(\xi;\xi_1^{})F_{\theta,N}^{}(\xi)~,
\end{array}
\label{jt33}
\end{equation}
where
\begin{equation}
\begin{array}{rcl}
&&F_{\theta_i^{},N_i^{}}^{}(\xi,\xi_1^{},
\ldots,\xi_{i-1}^{};\xi_i^{})\;:=\;(w_{\xi_{i-1}}^{}\otimes
w_{\xi_{i-1}}^{})\cdots (w_{\xi_{1}}^{}\otimes
w_{\xi_{1}}^{})(w_{\xi}^{}\otimes w_{\xi}^{})
F_{\theta_i^{},N_i^{}}^{}(\xi_i^{})
\\[9pt]
&&\phantom{aaaaaaaaaaaaa} \times (w_{\xi}^{-1}\otimes
w_{\xi}^{-1})(w_{\xi_{1}}^{-1} \otimes w_{\xi_{1}^{}}^{-1})\cdots
(w_{\xi_{i-1}^{}}^{-1}\otimes
w_{\xi_{i-1}^{}}^{-1})\quad(i=1,\ldots,k).
\end{array}
\label{jt34}
\end{equation}

\end{document}